\documentclass[12pt]{article}
\usepackage{amssymb}

\newtheorem{thm}     {Theorem}[section]
\newtheorem{prop}    [thm]{Proposition}

\newtheorem{cor}     [thm]{Corollary}
\newtheorem{lemma}   [thm]{Lemma}

\newcommand{\proof} {\noindent{\bf Proof. }}

\newcommand{\B}{\mathbb B}
\newcommand{\C}{\mathbb C}
\newcommand{\D}{\mathbb D}
\newcommand{\R}{\mathbb R}

\newcommand{\st}{{\rm st}}
\def\Re{{\rm Re\,}}
\def\Im{{\rm Im\,}}
\def\bar{\overline}
\def\area{{\rm Area}}

\begin{document}

\title{Gromov's Non-Squeezing Theorem
and Beltrami type equation}
\author{Alexandre Sukhov{*} and Alexander Tumanov{**}}
\date{}
\maketitle

{\small

* Universit\'e des Sciences et Technologies de Lille, Laboratoire
Paul Painlev\'e,
U.F.R. de
Math\'e-matique, 59655 Villeneuve d'Ascq, Cedex, France, sukhov@math.univ-lille1.fr

** University of Illinois, Department of Mathematics
1409 West Green Street, Urbana, IL 61801, USA, tumanov@illinois.edu
}
\bigskip

Abstract. We introduce a method for constructing J-complex discs. The method only uses the standard scheme for solving the Beltrami equation and the Schauder principle. As an application, we give a short self-contained proof of  Gromov's Non-Squeezing Theorem.
\bigskip

MSC: 32H02, 53C15.

Key words: almost complex structure, symplectic embedding, $J$-complex disc.

\section{Introduction}

We introduce a method for constructing $J$-complex
(pseudoholomorphic) discs for almost complex structures tamed by the standard symplectic form of $\C^n$.
The method only uses the standard scheme for solving the Beltrami equation (see \cite{Ve}) and the Schauder principle.
We do not need general machinery of pseudoholomorphic curves, in particular, compactness theorems, moduli spaces, etc. As an application we give a short proof of seminal Gromov's Non-Squeezing Theorem \cite{Gr}.

The proof of Gromov's theorem reduces to constructing
a proper $J$-complex disc of minimum area in the cylinder
$\D\times\C^{n-1}\subset \C^n$. Here $\D$ stands for the
unit disc in $\C$, and $J$ for a suitable almost
complex structure. This task in turn reduces to a boundary
value problem for an elliptic quasilinear system of PDE,
which is a vector analogue of the classical Beltrami equation.
The theory of such systems is well developed, especially
in the scalar case. However, most if not all general results
on the matter require linear boundary conditions whereas
in our problem, the condition that the first component
takes boundary values in the unit circle, is non-linear.
Our main idea was to replace the circular cylinder by the
triangular one, which does not matter in the original question.
Then the boundary conditions for the
sought $J$-complex disc become linear although with
discontinuous coefficients. The latter can be handled
by means of modified Cauchy-Green operators as we learned
from \cite{AM}, which inspired our work.
We hope this method will find other applications.

\section{Notation and terminology}
An almost complex structure $J$ on a smooth real manifold $M$, $\dim M = 2n$, is a map  which associates to every point $p \in M$ a linear isomorphism $J(p): T_pM \to T_pM$ of the tangent space $T_pM$ satisfying $J(p)^2 = -I$; here  $I$ denotes the identity map. A couple $(M,J)$ is called an almost complex manifold of complex dimension n.

Let $(M,J)$ and $(M',J')$ be almost complex manifolds. A $C^1$-map $f:M' \to M$ is called  $(J',J)$-complex or  $(J',J)$-holomorphic  if it satisfies {\it the Cauchy-Riemann equations}
\begin{eqnarray}
\label{CRglobal}
df \circ J' = J \circ df.
\end{eqnarray}
Denote by  $\D$  the unit disc in $\C$ and by $J_{\st}$
the standard complex structure
of $\C^n$; the value of $n$ will be clear from the context. For $M' = \D$ and $J' = J_{st}$,  we call a  map $f$ a $J$-{\it complex  disc}
(or a pseudoholomorphic disc).

Let $M$ be a smooth manifold of dimension $2n$. A closed non-degenerate exterior 2-form $\omega$ on $M$ is called a {\it symplectic form} on $M$. A pair $(M,\omega)$ is called a symplectic manifold. A basic example is $M = \C^{n}$ with the coordinates $Z_j = x_j + iy_j$, $j=1,...,n$. The form $\omega_{st} = \sum_{j=1}^n dx_j \wedge dy_j = \frac{i}{2}\sum_{j=1}^n dZ_j \wedge d\bar Z_j$ is called the standard symplectic form on $\C^n$.

A symplectic form $\omega$ {\it tames} an almost complex structure $J$ on $M$ if $\omega(u,Ju) > 0, \forall u \neq 0$.
A model example is provided by the standard symplectic form $\omega_{st}$ and the standard complex structure $J_{st}$ of $\C^n$. In this paper we deal only with the standard symplectic form so we denote it by $\omega$.

Let $J$ be an almost complex structure tamed by $\omega$ on $\C^n$. The Cauchy-Riemann equations (\ref{CRglobal}) for a $J$-complex disc $Z:\D\to\C^n$, $Z: \D \ni \zeta \mapsto Z(\zeta)$ can be rewritten in  the form
\begin{eqnarray}
\label{holomorphy}
Z_{\bar\zeta}=A(Z)\bar Z_{\bar\zeta},\quad
\zeta\in\D,
\end{eqnarray}
where a complex $n\times n$ matrix function $A=A(Z)$ satisfies the condition
\begin{eqnarray}
\label{taming}
\| A(Z) \| < 1 , \quad \forall Z \in \C^n.
\end{eqnarray}
Here the matrix norm is induced by the Euclidean inner product and $J$ is tamed by $\omega$ if and only if (\ref{taming}) holds. In fact, $A$ is uniquely determined by $J$ as the matrix representation of  the complex anti-linear operator $(J_\st+J)^{-1}(J_\st-J)$. In particular, $A(Z) = 0$ if $J(Z) = J_{st}$. Conversely, every $A$ satisfying (\ref{taming}) defines a unique almost complex structure tamed by $\omega$, see \cite{Aud}. We call $A$ the complex matrix of $J$. Thus, $J$-complex discs are precisely solutions of the system (\ref{holomorphy}), which is elliptic by (\ref{taming}).

The system (\ref{holomorphy}) generalizes the
classical Beltrami equation (see \cite{Ve})
to higher dimension. It still makes sense for
$Z \in W^{1,p}(\D)$ for $p > 2$.
Here $W^{1,p}(\D)$, $p \in [1,\infty]$, stands for the Sobolev class of maps $\D\to\C^m$ with first partial derivatives in $L^p$; the value of $m$ will be clear from context.
For $p>2$ the space $W^{1,p}(\D)$ is continuously embedded into $C^\alpha(\overline\D)$, the space of maps satisfying the H\"older condition with exponent $\alpha=1-2/p$, in particular,
elements of $W^{1,p}(\D)$ are continuous in $\bar\D$.

For a map $Z: \D \to \C^n$, the (symplectic) {\it area} of $Z$ is given by
\begin{eqnarray}
\label{area}
\area(Z) = \int_\D Z^*\omega.
\end{eqnarray}
In the case where $Z$ is a $J$-complex disc, it coincides with the area induced by the Riemannian metric canonically defined by $J$ and $\omega$; in particular, it coincides with the Euclidean area if $J = J_{st}$ (see, for instance, \cite{Aud}). We use the same notation for the Euclidean area of complex analytic sets
in $\C^n$.

\section{Results}

Denote by $\Delta$ the triangle $\Delta = \{ z \in \C: 0 < \Im z < 1 - \vert \Re z \vert \}$.
Note that   $\area(\Delta) = 1$.
Consider the cylinder $\Sigma = \Delta \times \C^{n-1}$ in $\C^n$.
We now use the notation
$$
Z = (z,w) = (z,w_1,...,w_{n-1})  \in \C \times \C^{n-1}=\C^n
$$
for the coordinates in $\C^n$.
Our main result is the following

\begin{thm}
\label{ThDiscs}
Let $A$ be a continuous $n \times n$ matrix function on $\C^n$ vanishing on $\C^n\setminus\Sigma$. Suppose there is a constant $0<a<1$ such that
\begin{eqnarray}
\label{norm}
\| A(Z) \|  \leq a, \quad  \forall Z \in \Sigma.
\end{eqnarray}
Then there exists $p > 2$ such that for every point $(z^0,w^0) \in \Sigma$ there is a solution $Z\in W^{1,p}(\D)$ of (\ref{holomorphy}) such that $Z(\bar\D)\subset\bar\Sigma$, $(z^0,w^0) \in Z(\D)$, $\area(Z) = 1$, and
\begin{eqnarray}
\label{BC0}
Z(b\D) \subset
b\Sigma = (b\Delta) \times \C^{n-1}.
\end{eqnarray}
\end{thm}

As a consequence we obtain Gromov's Non-Squeezing Theorem.
We essentially repeat Gromov's \cite{Gr} argument that consists of constructing a $J$-complex curve of small area and pulling it back. However, we use $J$-complex discs instead of compact curves.
Denote by $\B^n$ the Euclidean unit ball in $\C^n$.

\begin{cor}
\label{squiz}
Let $G$ be a domain  in $R\D \times \C^{n-1}$ where $R > 0$. Suppose that $r> 0$ and there exists a $C^1$-diffeomorphism $\Phi: r\B^n \to G$ with $\Phi^* \omega = \omega$. Then $r \leq R$.
\end{cor}
\proof A diffeomorphism whose $z$-component is an area-preserving map and whose $w$-components are the identity maps, preserves the form $\omega$. This observation reduces the proof to the case where $G$ is contained in the cylinder
$\Sigma_R := \sqrt{\pi}R\Delta \times \C^{n-1}$.
Since $\Phi^*\omega=\omega$, then the almost complex structure $J := d\Phi \circ J_{st} \circ d\Phi^{-1}$ is tamed by $\omega$.
Then the complex matrix $\tilde A$ of $J$ satisfies
$||\tilde A(Z)||<1$ for $Z\in G$.
Fix $\epsilon > 0$.
Let $\chi$ be a smooth cut-off function with support in $G$
and such that $\chi=1$ on $\Phi((r-\epsilon)\bar\B^n)$.
Define $A=\chi\tilde A$. Let $p = \Phi(0)$.
Since $J$ is continuous in $G$, then there is a constant
$a<1$ such that (\ref{norm}) holds for $A$.
By Theorem \ref{ThDiscs} there exists a solution $Z$ of (\ref{holomorphy}) such that $p \in Z(\D)$, $Z(b\D) \subset b\Sigma_R$ and $\area(Z) = \pi R^2$.
Then $X = \Phi^{-1}(Z(\D)) \cap (r - \epsilon)\B^n$ is a closed $J_{st}$-complex curve  in $(r-\epsilon)\B^n$. Furthermore, $0 \in X$ and $\area(X) \leq \pi R^2$.
On the other hand, by the classical result due to Lelong (see, e.g., \cite{Ch}) we have $\area(X) \geq \pi (r-\epsilon)^2$. Since $\epsilon$ is arbitrary, then $r\le R$ as desired. $\blacksquare$
\medskip

In the rest of the paper we prove Theorem \ref{ThDiscs}.

\section{Modified Cauchy-Green operators}

Introduce the functions
$$
R(\zeta) = e^{3\pi i/4}(\zeta - 1)^{1/4} (\zeta + 1)^{1/4}(\zeta - i)^{1/2}, \qquad
X(\zeta)= R(\zeta)/\sqrt{\zeta}.
$$
Here we choose the branch of $R$ continuous in $\bar\D$
with $R(0) = e^{3\pi i/4}$.
We need the function $X$ only on the circle $b\D$.
We do not care about the sign of $X$, nevertheless,
for definiteness, we choose the branch of $\sqrt{\zeta}$
continuous in $\C$ with deleted positive real line,
$\sqrt{-1}=i$.
One can see that $(X(\zeta))^4\in\R$ for $\zeta\in b\D$.
Then $\arg X$ is constant on each arc
$\gamma_1 = \{ e^{i\theta} : 0 < \theta < \pi/2 \}$, $\gamma_2 = \{ e^{i\theta} : \pi/2 < \theta < \pi \}$, $\gamma_3 = \{ e^{i\theta} : \pi < \theta < 2\pi\}$.
Moreover, $\arg X$ on these arcs is equal to  $3\pi/4$, $\pi/4$ and $0$ respectively.  Therefore, the function $X$ satisfies the boundary conditions
\begin{eqnarray}
\label{BC}
 \left\{ \begin{array}{cccc}
& &\Im (1+i)X(\zeta) = 0, \quad \zeta \in \gamma_1,\\
& &\Im (1-i)X(\zeta) = 0, \quad \zeta \in \gamma_2,\\
& &\Im X(\zeta) = 0, \quad \zeta \in \gamma_3,
\end{array} \right.
\end{eqnarray}
which represent the lines through 0 parallel to the sides of
the triangle $\Delta$.

We will use modifications of the classical Cauchy-Green operator
$$T(f)(\zeta) = \frac{1}{2\pi i} \int_\D \frac{f(t) dt \wedge d\overline{t}}{t-\zeta}.$$
Recall that $T: L^p(\D) \to W^{1,p}(\D)$ is bounded for $p > 2$ and $(\partial/\partial\overline{\zeta}) Tf = f$ as Sobolev's derivative, i.e., $T$ solves the $\overline\partial$-problem in $\D$. Furthermore, $Tf$ is holomorphic on $\C \setminus \overline{\D}$, see \cite{Ve}.

Let $Q$ be a function in $\D$. We call it a weight function. Introduce the operator
\begin{eqnarray*}
& &T_Qf(\zeta) = Q(\zeta)\left ( T(f/Q)(\zeta) + \zeta^{-1}\overline{T(f/Q)(1/\overline{\zeta})}\right )\\
& &= Q(\zeta) \int_\D \left ( \frac{f(t)}{Q(t)(t-\zeta)} + \frac{\overline{f(t)}}{\overline{Q(t)}(\overline{t}\zeta - 1)} \right )\frac{dt \wedge d\bar t}{2\pi i}.
\end{eqnarray*}
We will need only the operators corresponding to two special weights, namely $T_1 = T_Q$ with $Q = \zeta - 1$ and $T_2 = T_Q$ with $Q = R$.
We also define formal derivatives
$S_jf(\zeta) = (\partial/\partial\zeta) T_jf(\zeta)$ as integrals in the sense of the Cauchy principal value.
The operator $T_1$ was first introduced by Vekua \cite{Ve}
whereas operators similar to $T_2$ apparently were
first introduced by Antoncev and Monakhov \cite{AM,Mo}
for application to problems of gas dynamics.
The operators $T_j$ and $S_j$, $j=1,2$, have the following
properties, see \cite{Mo, Ve}.
\begin{prop}
\label{OpBounVal}
\begin{itemize}
\item[(i)] Each $S_j :L^p(\D) \to L^p(\D)$, $j=1,2$, is a bounded linear operator for $p_1 < p < p_2$. Here for $S_1$ one has $p_1 = 1$ and $p_2 = \infty$ and for $S_2$ one has $p_1 = 4/3$ and $p_2 = 8/3$. For $2<p<p_2$, one has $S_jf(\zeta) = (\partial/\partial\zeta) T_jf(\zeta)$ as Sobolev's derivatives.
\item[(ii)]  Each $T_j :L^p(\D) \to W^{1,p}(\D)$, $j=1,2$,  is a bounded linear operator for $2<p<p_2$. In particular, $T_j:L^p(\D) \to L^\infty(\D)$ is a compact operator. For $f \in L^p(\D)$, $2<p<p_2$, one has $(\partial/\partial\overline{\zeta}) T_j f = f$ on $\D$ as Sobolev's derivative.
\item[(iii)]  For every $f \in L^p(\D)$, $p > 2$, the function $T_1f$ satisfies $\Re T_1f\vert_{b\D} = 0$
    whereas $T_2f$ satisfies
    the same boundary conditions (\ref{BC}) as $X$.
\item[(iv)] Each $S_j: L^2(\D) \to L^2(\D)$, $j=1,2$, is an isometry.
\item[(v)]  The function $p \mapsto \| S_j \|_{L^p}$ approaches $\| S_j \|_{L^2} = 1$ as $p\searrow 2$.
\end{itemize}
\end{prop}
\proof
(i, ii) See \cite{Mo, Ve}.

(iii) If $|\zeta|=1$, then
$$T_2f(\zeta) = (R(\zeta)/\sqrt\zeta)\left (\sqrt{\zeta}T(f/R)(\zeta) + \overline{\sqrt{\zeta}T(f/R)(\zeta)}\right ).$$
Since the expression in parentheses is real, then
$\arg T_2f(\zeta) = \arg (\pm X(\zeta))$, hence
the conclusion. The proof for $T_1$ is similar,
but simpler.

(iv) This is proved in \cite{Mo} in a more general situation.
For completeness, we give a simple proof in our special cases. Let $f$ be a smooth function with compact support in $\D$. Since $T_jf(b\D)$ lies on finitely many lines, then $\area(T_jf)=0$. Therefore by Stokes' formula
\begin{eqnarray*}
& &0 = (i/2)\int_{b\D} T_j f d\overline{T_jf} = (i/2)\int_{\D}dT_jf \wedge d\overline{T_jf}\\
& &= (i/2)\int_\D (S_j f d\zeta + f d \overline\zeta) \wedge (\overline{S_jf}d\overline{\zeta} + \overline{f}d\zeta) \\
& &= (i/2)\int_\D \vert S_j f \vert^2 d\zeta \wedge d \overline{\zeta} - (i/2)\int_\D \vert f \vert^2 d\zeta \wedge d\overline\zeta.
\end{eqnarray*}
Hence $\| S_jf \|_{L^2(\D)} = \| f \|_{L^2(\D)}$ and by density this equality holds for all $f \in L^2(\D)$.

(v) This follows by  the Riesz-Thorin interpolation theorem. $\blacksquare$

\section{Reduction to an integral equation}

Consider the biholomorphism $\Phi:\D \to \Delta$ satisfying $\Phi(\pm 1) = \pm 1$ and $\Phi(i) = i$. Note that $\Phi \in W^{1,p}(\D)$ for $p> 2$ close enough to $2$ by the classical results on boundary behavior of conformal maps.
Following the standard scheme for solving the Beltrami
equation \cite{Ve},
we look for a solution $Z = (z,w):\D \to \C^{n}$  of (\ref{holomorphy}) in the form
\begin{eqnarray}
\label{MR}
\left\{ \begin{array}{cccc}
& & z = T_2u + \Phi,\\
& & w = T_1v - T_1v(\tau) + w^0.
\end{array} \right.
\end{eqnarray}
for some $\tau \in \D$; hence, $w(\tau) = w^0$.
The Cauchy-Riemann equation (\ref{holomorphy}) for $Z$
of the form (\ref{MR}) turns into the integral equation
\begin{eqnarray}
\label{mainsystem}
\left(
\begin{array}{cl}
 u\\
 v
\end{array}
\right) = A(z,w)\left(
\begin{array}{cl}
 \overline{S_2u} +\overline{\Phi'}\\
 \overline{S_1v}
\end{array}
\right).
\end{eqnarray}
Our task reduces to showing that there exists a solution of (\ref{MR}, \ref{mainsystem}) so that $z(\tau)=z^0$ for some $\tau\in\D$.
We first obtain a priori estimates for (\ref{MR}, \ref{mainsystem}). After increasing the constant $a<1$ in (\ref{norm}) if necessary, we assume that for all $p$ close to 2
and all $Z\in\C^n$ we have $||A(Z)||_p\le a<1$. Here
$||\,.\,||_p$ stands for the matrix norm induced by the $p$-norm
in $\C^n$.

Using that $s = \max_j \| S_j \|_{L^p}\to 1$ as $p\searrow 2$,
we choose $p > 2$ close to $2$ such that $as<1$. Then for every fixed $Z = (z,w):\D \to \C^{n}$, by the contraction principle
in $L^p(\D)$, the equation (\ref{mainsystem}) has a unique solution $Y = (u,v)$ satisfying
$$
\| Y \|_{L^p} \leq a \left (s\| Y \|_{L^p} + \| \Phi \|_{L^p}\right ), \qquad
\| Y \|_{L^p} \leq M_1:=\frac{a\| \Phi' \|_{L^p}}{1 - a s}.
$$
It follows by (\ref{MR}) that there exists a constant $M > 0$ depending on $M_1$ and $w^0$ such that
\begin{eqnarray}
\label{ball1}
\| z \|_{L^\infty} \leq M, \quad \| w \|_{L^\infty} \leq M.
\end{eqnarray}
We now define a continuous map $\Psi:\C \to \overline{\D}$
\begin{eqnarray*}
\Psi(z) = \left\{ \begin{array}{cccc}
& & \Phi^{-1}(z), \quad z \in \overline{\Delta},\\
& & \Phi^{-1}(b\Delta \cap [z^0,z]), \quad  z \in \C \setminus  \overline{\Delta}.
\end{array} \right.
\end{eqnarray*}
Here $[z^0,z]$ is the line segment from $z^0$ to $z$, and the intersection $b\Delta \cap [z^0,z]$ consists of a single point. (Note that in the definition of $\Psi$ one can replace the point $z^0$ by a fixed point, say $i/2$, making the function $\Psi$ independent of the initial data $z^0$. We use the point $z^0$ for convenience of presentation.)

Consider the balls $E_z = \{ z \in L^\infty(\D): \| z \|_{L^\infty} \leq M \}$
and $E_w = \{ w \in L^\infty(\D): \| w \|_{L^\infty} \leq M \}$
and define $E = E_z \times E_w \times \overline{\D}$. Introduce the map $F: E \to E$,
$F:(z,w,\tau) \mapsto (\tilde{z},\tilde{w},\tilde{\tau})$
defined by
\begin{eqnarray*}
& &\tilde{z} = T_2u+\Phi,\\
& &\tilde{w} = T_1v - T_1v(\tau) + w_0,\\
& &\tilde{\tau} = \Psi(z^0 - T_2u(\tau)).
\end{eqnarray*}
Here $(u,v)$ is a solution of (\ref{mainsystem}).
The map $F$ is continuous because $A$ is.
The set $E$ is convex and the operators $T_j: L^{p}(\D) \to L^\infty(\D)$ are compact. It follows now by Schauder's principle that the map $F$ has a fixed point $(z,w,\tau)$. It satisfies (\ref{MR}), (\ref{mainsystem}) and $\tau = \Psi(z^0 - T_2u(\tau))$.

\section{Properties of the solution}

By (\ref{MR}) and (\ref{mainsystem}), the map $Z=(z,w)\in W^{1,p}(\D)$ satisfies the Cauchy-Riemann equations (\ref{holomorphy}) and $w(\tau)=w^0$. We now prove the other conclusions
of Theorem \ref{ThDiscs}.

\begin{lemma}
\label{tau-lemma}
$\tau \in \D$ and $z(\tau) = z^0$.
\end{lemma}
\proof
Suppose otherwise that $\tau \in b\D$.
Then $z^0 -T_2u(\tau)\notin\Delta$, in particular,
$q:= T_2u(\tau) = z(\tau) - \Phi(\tau)\ne 0$.
By the definition of the map $\Psi$, we have
$\Phi(\tau) = b\Delta \cap [z^0,z^0-q]$.
For definiteness, suppose $\tau\in\bar\gamma_1$.
Then $\Phi(\tau) \in [1,i]$, and by the boundary conditions
(\ref{BC}) for $T_2$, $q$ is a real multiple of $1-i$.
Now we must have
$\Phi(\tau) = [1,i] \cap [z^0,z^0-q]$,
which is absurd because the lines $[1,i]$ and $[z^0,z^0-q]$
are parallel.

Now since $\tau \in \D$, then $\tau = \Psi(z^0 - T_2u(\tau))$ implies $\Phi(\tau) = z^0 - T_2u(\tau)$, and $z(\tau) = T_2u(\tau) + \Phi(\tau) = z^0$.
$\blacksquare$

\begin{lemma}
\label{image}
The map $z$ satisfies $z(\bar\D) \subset \bar\Delta$,
$z(b\D) \subset b\Delta$, and $\deg z=1$;
here $\deg z$ denotes the degree of the map
$z|_{b\D}:b\D \to b\Delta$.
In particular, $Z$ satisfies (\ref{BC0}).
\end{lemma}
\proof
Let $G=\{\zeta\in\D: z(\zeta)\notin\bar\Delta\}$.
Arguing by contradiction, suppose $G\ne\emptyset$.
Since $z$ is continuous, then $G$ is open.
Let $G_1$ be a non-empty connected component of $G$.
Then $z(bG_1)\subset\bar\Delta$.
Since $A=0$ on $\C^n\setminus\Sigma$ and by
(\ref{holomorphy}), $z$ is holomorphic on $G_1$.
But then the set $z(G_1)$ has the farthest point
from $\bar\Delta$, which violates the maximum principle.
Hence $G=\emptyset$, $z(\D) \subset \bar\Delta$,
and by continuity $z(\bar\D) \subset \bar\Delta$.

By the boundary properties of $T_2u$ and $\Phi$,
the map $z=T_2u+\Phi$ takes the arcs $\gamma_j$, $j=1,2,3$,
to the lines containing the corresponding sides of the triangle $\Delta$. Since $z(\bar\D) \subset \bar\Delta$, then
the images $z(\gamma_j)$, $j=1,2,3$, are exactly the sides
of $\Delta$. Hence $z(b\D) \subset b\Delta$
and $\deg z=1$. $\blacksquare$

\begin{lemma}
\label{area1}
$\area(Z)=1$.
\end{lemma}
\proof
By Stokes' formula
$$
\area(Z)
=\frac{i}{2} \,\int_\D
\left(dz\wedge d\bar z
+\sum dw_j\wedge d\bar{w_j}\right)
=\frac{i}{2}\,
\int_{b\D} z\, d\bar z
+\sum \frac{i}{2}\,\int_{b\D} w_j\, d\bar{w_j}.
$$
We now evaluate each integral on the right separately.
Since $z(b\D) \subset b\Delta$, and $\deg z=1$, then
$(i/2)\int_{b\D} z\, d\bar z = \area(\Delta)=1$.
By the boundary properties of $T_1$, the real part
$\Re w_j=\Re w_j^0$ is constant, therefore
$\int_{b\D} w_j\, d\bar{w_j}=0$.
Hence $\area(Z)=1$ as desired.

We point out that although $Z$ has fairly low regularity,
the use of Stokes' formula is legitimate.
Indeed, we can approximate $(u,v)$ in $L^p(\D)$ by smooth functions with compact support in $\D$ and define the approximation $\tilde Z$ of $Z=(z,w)$ by (\ref{MR}).
By the boundary properties of $T_1$ and $T_2$ and
the above argument, $\area(\tilde Z)=1$. Since
$\tilde Z$ approaches $Z$ in $W^{1,p}(\D)$,
then $\area(Z)=1$.
$\blacksquare$
\medskip

The proof of Theorem \ref{ThDiscs} is complete.

\section{Why triangle?}

The reader may wonder why we choose a triangle as
the base of the cylinder in Theorem \ref{ThDiscs}.
We add a few lines on this matter.

First of all, we restrict to convex polygons because
the proof of Lemma \ref{tau-lemma} needs convexity.
For a convex polygon other than a triangle
the construction of the operator $T_2$ described
in Section 4 does not go through.
Indeed, consider, say a quadrilateral $K$ with angles
$\pi \alpha_j$, $1\le j \le 4$.
Then following the construction, we put
$R(\zeta)=\sigma \prod (\zeta-\zeta_j)^{\alpha_j}$,
$\zeta_j\in b\D$. Then $\sum \alpha_j=2$,
and for a suitable constant $\sigma\ne0$, the function $X(\zeta)=R(\zeta)/\zeta$ satisfies the desired boundary conditions. Then we define
$T_2 f(\zeta) = R(\zeta)\left ( T(f/R)(\zeta) + \zeta^{-2}\overline{T(f/R)(1/\overline{\zeta})}\right )$,
which satisfies the same boundary conditions as $X$.
However, $T_2 f$ defined that way is not even in
$L^p$ for any $p>2$.

Finally, for the unit square
$K=\{z\in \C: |\Re z|<1, |\Im z|<1 \}$
the analogue of $T_2$ clearly does not exist.
Otherwise we take a conformal map $\Phi:\D\to K$
and define
$f=-T_2 \bar {\Phi'} + \bar\Phi$.
Then $f:\D\to K$ is holomorphic and continuous up to
the boundary, but the degree of the map
$f|_{b\D}:b\D\to bK$ is negative, which is absurd.

\end{document}